\documentclass[12pt]{article}
\usepackage{latexsym}
\usepackage{graphicx}
\usepackage{cite}
\usepackage{enumitem}
\usepackage{amssymb,mathrsfs}
\usepackage{booktabs}
\usepackage[justification=centering]{caption}
\usepackage{subfigure}
\usepackage{float}

\newtheorem{theorem}{Theorem}[part]

\newtheorem{lemma}{Lemma}[part]
\newtheorem{corollary}{Corollary}[part]

\def \ep{\hbox{ }\hfill$\Box$}

\begin{document}
\title{Stationary probability vectors of higher-order two-dimensional transition probability tensors
}

\author{Zheng-Hai Huang\footnote{%
    School of Mathematics, Tianjin University, Tianjin 300354, P.R. China ({\tt huangzhenghai@tju.edu.cn}).
    This author was supported by the National Natural Science Foundation of China (Grant No. 11431002).}
  \and Liqun Qi\footnote{%
    Department of Applied Mathematics, The Hong Kong Polytechnic University,
    Hung Hom, Kowloon, Hong Kong ({\tt maqilq@polyu.edu.hk}).
    This author's work was partially supported by the Hong Kong Research Grant Council
    (Grant No. PolyU  15302114, 15300715, 15301716 and 15300717).}
    }
\date{February 28, 2018}

\maketitle

\begin{abstract}
\noindent
In this paper we investigate stationary probability vectors of higher-order two-dimensional symmetric transition probability tensors. We show that there are two special symmetric transition probability tensors of order $m$ dimension 2, which have and only have two stationary probability vectors; and any other symmetric transition probability tensor of order $m$ dimension 2 has a unique stationary probability vector. As a byproduct, we obtain that any symmetric transition probability tensor of order $m$ dimension 2 has a unique positive stationary probability vector; and that any symmetric irreducible transition probability tensor of order $m$ dimension 2 has a unique stationary probability vector.
\vspace{3mm}

\noindent {\bf Key words:}\hspace{2mm} Transition probability tensor; higher-order Markov chain; stationary probability vector; eigenvalue of tensor. \vspace{3mm}

\noindent {\bf Mathematics Subject Classifications(2000):}\hspace{2mm} 15A18; 15A69; 65F15; 60J10; 60J22. \vspace{3mm}

\end{abstract}

\newpage

\section{Introduction}
\setcounter{equation}{0} \setcounter{assumption}{0}
\setcounter{theorem}{0} \setcounter{proposition}{0}
\setcounter{corollary}{0} \setcounter{lemma}{0}
\setcounter{definition}{0} \setcounter{remark}{0}
\setcounter{algorithm}{0}

It is well known that higher-order Markov chains have various applications in many areas \cite{BP94,MZ97,Raf85,RT94,Wat95,CNg06}. An $(m-1)$-order $n$-dimensional Markov chain is basically characterized by its associated nonnegative tensor $\mathscr{P}$ which is an $m$-order $n$-dimensional tensor with entries $p_{i_1i_2\cdots i_m}\geq 0$ for all $i_j\in \{1,2,\ldots,n\}$ and $j\in \{1,2,\ldots,m\}$ satisfying
\begin{eqnarray}\label{e-sect1-1}
0\leq p_{i_1i_2\cdots i_m}=\mbox{\rm Prob}(X_{t+1}=i_1|X_t=i_2,\ldots,X_{t-m+2}=i_m)\leq 1
\end{eqnarray}
where $\{X_t: t=0,1,\ldots\}$ represents the stachastic process that takes on $n$ states $\{1,2,\ldots,n\}$, and for any $i_2,\ldots,i_m\in \{1,2,\ldots,n\}$,
\begin{eqnarray}\label{e-sect1-2}
\sum_{i_1=1}^np_{i_1i_2\cdots i_m}= 1.
\end{eqnarray}

We will use $\mathbb{T}_{m,n}$ to denote the set of all $m$-order $n$-dimensional real tensors. For any $\mathscr{P}=(p_{i_1\cdots i_m})\in \mathbb{T}_{m,n}$, if the entries $p_{i_1i_2 \ldots i_m}$ are invariant under any permutation of their indices, then $\mathscr{P}$ is called a {\bf symmetric tensor}.
A tensor $\mathscr{P}\in \mathbb{T}_{m,n}$ is called a {\bf transition probability tensor} if it satisfies (\ref{e-sect1-1}) and (\ref{e-sect1-2}). A vector
$$
x^*\in \left\{x\in\mathbb{R}^n: x_i\geq 0\;\mbox{\rm for all}\; i\in \{1,2,\ldots,n\}\; \mbox{\rm and}\; \sum_{i=1}^nx_i=1\right\}
$$
is called a {\bf stationary probability vector} of  $\mathscr{P}\in \mathbb{T}_{m,n}$ if
$$
\sum_{i_2,\ldots,i_m=1}^np_{ii_2\cdots i_m}x_{i_2}^*\cdots x_{i_m}^*=x_{i}^*
$$
holds for all $i\in \{1,2,\ldots,n\}$, which is just a $Z_1$-eigenvector associated with the $Z_1$-eigenvalue 1 \cite{CZ13}. It is also closely related to the $Z$-eigenvector of a tensor \cite{Qi05,Lim05}.

Transition probability tensors and the associated stationary probability vectors are important issues in studies of higher-order Markov chains \cite{NgQZ09,CZ13,HQ14,LNg14,LZ16,BB16,CPZ17}. In particular, the uniqueness of the stationary probability vector of transition probability tensors has attracted a lot of interest. Li and Ng \cite{LNg14} proposed some conditions which ensure the uniqueness of the stationary probability vector of transition probability tensors and the linear convergence of the proposed iterative method; Hu and Qi \cite{HQ14} studied the uniqueness of the stationary probability vector of the third order $n$-dimensional positive transition probability tensor, and they proved that an irreducible transition probability tensor of order 3 dimension 2 has a unique stationary probability vector; and Chang and Zhang \cite{CZ13} investigated sufficient
conditions for transition probability tensors to ensure the uniqueness of the stationary probability vector by using three different methods: contraction mappings, monotone operators, and the Brouwer index of fixed points.

More recently, Culp, Pearson and Zhang \cite{CPZ17} investigated symmetric irreducible transition probability tensors of order 4 dimension 2 and of order 3 dimension 3, and showed that a symmetric irreducible transition probability tensor in these orders and dimensions has a unique stationary probability vector.

In this paper, we give a full characterization on the stationary probability vectors of $m$-order $2$-dimensional symmetric transition probability tensors. We show that there are two special symmetric transition probability tensors of order $m$ dimension 2, which have and only have two stationary probability vectors, where one is $(\frac 12,\frac 12)^\top$, and the other one is $(1,0)^\top$ or $(0,1)^\top$. In particular, we show that if the concerned symmetric transition probability tensor of order $m$ dimension 2 is not one of the above two tensors, then it has a unique stationary probability vector, which is $(\frac 12,\frac 12)^\top$.
As a byproduct, we obtain that any symmetric irreducible transition probability tensor of order $m$ dimension 2 has a unique stationary probability vector. When $m=4$, such a result was obtained by Culp, Pearson and Zhang in \cite{CPZ17}.

Throughout this paper, we assume that $m\geq 3$ is an integer number.

\section{Main results}
\setcounter{equation}{0} \setcounter{assumption}{0}
\setcounter{theorem}{0} \setcounter{proposition}{0}
\setcounter{corollary}{0} \setcounter{lemma}{0}
\setcounter{definition}{0} \setcounter{remark}{0}
\setcounter{algorithm}{0}\setcounter{example}{0}

Let $\mathscr{P}=(p_{i_1i_2\cdots i_m})\in \mathbb{T}_{m,2}$ be a transition probability tensor. Then $z=(z_1,z_2)^\top\in \mathbb{R}^2$ is a stationary probability vector of $\mathscr{P}$ if and only if $z_1,z_2\geq 0$, $z_1+z_2=1$, and
\begin{eqnarray}\label{e-func-f1}
f_1(z_1,z_2):=\sum\limits_{i_2,\ldots,i_m=1}^2p_{1i_2\cdots i_m}z_{i_2}\cdots z_{i_m}=z_1,
\end{eqnarray}\vspace{-3mm}
\begin{eqnarray}\label{e-func-f2}
f_2(z_1,z_2):=\sum\limits_{i_2,\cdots,i_m=1}^2p_{2i_2\cdots i_m}z_{i_2}\cdots z_{i_m}=z_2.
\end{eqnarray}
In the following, we denote $x:=z_1$ and $y:=z_2$, and
\begin{eqnarray}\label{e-Delta}
\Delta:=\left\{(x,y)^\top\in \mathbb{R}^2: x\geq 0, y\geq 0, x+y=1\right\}.
\end{eqnarray}
Then, $(x,y)^\top\in \mathbb{R}^2$ is a stationary probability vector of $\mathscr{P}$ if and only if $(x,y)^\top\in \Delta$ and
\begin{eqnarray}\label{e-system-1}
f_1(x,y)=x\quad\mbox{\rm and}\quad f_2(x,y)=y.
\end{eqnarray}
Throughout this paper, we denote $a:=p_{1\cdots 11}$ and $b:=p_{21\cdots 11}$.

\begin{lemma}\label{lem-1}
Suppose that $\mathscr{P}=(p_{i_1i_2\cdots i_m})\in \mathbb{T}_{m,2}$ is a symmetric transition probability tensor. Let $f_1(\cdot,\cdot)$, $f_2(\cdot,\cdot)$, and $\Delta$ be defined by (\ref{e-func-f1}), (\ref{e-func-f2}), and (\ref{e-Delta}), respectively. Then, for any $(x,y)^\top\in \Delta$, we have the following results.
\begin{itemize}
\item[(i)] If $m$ is an even number, then
\begin{eqnarray*}
\begin{array}{rcl}
f_1(x,y)&=&ax^{m-1}+bC_{m-1}^1x^{m-2}y+aC_{m-1}^2x^{m-3}y^2+bC_{m-1}^3x^{m-4}y^3+\cdots+\\
&& ax^3y^{m-4}+bC_{m-1}^{m-3}x^2y^{m-3}+aC_{m-1}^{m-2}xy^{m-2}+by^{m-1},\vspace{2mm}\\
f_2(x,y)&=&bx^{m-1}+aC_{m-1}^1x^{m-2}y+bC_{m-1}^2x^{m-3}y^2+aC_{m-1}^3x^{m-4}y^3+\cdots+\\
&& bx^3y^{m-4}+aC_{m-1}^{m-3}x^2y^{m-3}+bC_{m-1}^{m-2}xy^{m-2}+ay^{m-1}.
\end{array}
\end{eqnarray*}
\item[(ii)] If $m$ is an odd number, then
\begin{eqnarray*}
\begin{array}{rcl}
f_1(x,y)&=&ax^{m-1}+bC_{m-1}^1x^{m-2}y+aC_{m-1}^2x^{m-3}y^2+bC_{m-1}^3x^{m-4}y^3+\cdots+\\
&& ax^4y^{m-5}+bx^3y^{m-4}+aC_{m-1}^{m-3}x^2y^{m-3}+bC_{m-1}^{m-2}xy^{m-2}+ay^{m-1},\vspace{2mm}\\
f_2(x,y)&=&bx^{m-1}+aC_{m-1}^1x^{m-2}y+bC_{m-1}^2x^{m-3}y^2+aC_{m-1}^3x^{m-4}y^3+\cdots+\\
&& bx^4y^{m-5}+ax^3y^{m-4}+bC_{m-1}^{m-3}x^2y^{m-3}+aC_{m-1}^{m-2}xy^{m-2}+by^{m-1}.
\end{array}
\end{eqnarray*}
\end{itemize}
\end{lemma}

\noindent{\bf Proof.} Since $\mathscr{P}$ is a symmetric tensor, the above equalities can be rewritten as
\begin{eqnarray}\label{e-lem-1-1}
\begin{array}{rcl}
f_1(x,y)&=&p_{1\cdots 11}x^{m-1}+p_{1\cdots 12}C_{m-1}^1x^{m-2}y+p_{1\cdots 122}C_{m-1}^2x^{m-3}y^2\\
&&+p_{1\cdots 1222}C_{m-1}^3x^{m-4}y^3+\cdots+p_{11112\cdots 2}x^3y^{m-4}\\
&&+p_{1112\cdots 2}C_{m-1}^{m-3}x^2y^{m-3}+p_{112\cdots 2}C_{m-1}^{m-2}xy^{m-2}+p_{12\cdots 2}y^{m-1},\vspace{2mm}\\
f_2(x,y)&=&p_{21\cdots 11}x^{m-1}+p_{21\cdots 12}C_{m-1}^1x^{m-2}y+p_{21\cdots 122}C_{m-1}^2x^{m-3}y^2\\
&&+p_{21\cdots 1222}C_{m-1}^3x^{m-4}y^3+\cdots+p_{21112\cdots 2}x^3y^{m-4}\\
&&+p_{2112\cdots 2}C_{m-1}^{m-3}x^2y^{m-3}+p_{212\cdots 2}C_{m-1}^{m-2}xy^{m-2}+p_{22\cdots 2}y^{m-1}=y.
\end{array}
\end{eqnarray}
Suppose that $m$ is an even number. Since $\mathscr{P}$ is a transition probability tensor, it follows that
\begin{eqnarray*}
\begin{array}{l}
p_{1\cdots 11}+p_{1\cdots 12}=1,\;p_{1\cdots 12}+p_{1\cdots 122}=1,\; p_{1\cdots 122}+p_{1\cdots 1222}=1,\;\ldots, \\
\qquad p_{11112\cdots 2}+p_{1112\cdots 2}=1,\;p_{1112\cdots 2}+p_{112\cdots 2}=1,\;
p_{112\cdots 2}+p_{12\cdots 2}=1,\vspace{2mm}\\
p_{21\cdots 11}+p_{21\cdots 12}=1,\;p_{21\cdots 12}+p_{21\cdots 122}=1,\;
p_{21\cdots 122}+p_{21\cdots 1222}=1,\;\ldots,\\
\qquad p_{21112\cdots 2}+p_{2112\cdots 2}=1,\;p_{2112\cdots 2}+p_{212\cdots 2}=1,\;
p_{212\cdots 2}+p_{22\cdots 2}=1,
\end{array}
\end{eqnarray*}
and hence,
\begin{eqnarray*}
\begin{array}{l}
p_{1\cdots 11}=p_{1\cdots 122}=\cdots=p_{11112\cdots 2}=p_{112\cdots 2}=a,\vspace{2mm}\\
p_{1\cdots 12}=p_{1\cdots 1222}=\ldots=p_{1112\cdots 2}=p_{12\cdots 2}=b,\vspace{2mm}\\
p_{21\cdots 11}=p_{21\cdots 122}=\cdots=p_{21112\cdots 2}=p_{212\cdots 2}=b,\vspace{2mm}\\
p_{21\cdots 12}=p_{21\cdots 1222}=\ldots=p_{2112\cdots 2}=p_{22\cdots 2}=a.
\end{array}
\end{eqnarray*}
Thus, two equalities given in (i) hold from (\ref{e-lem-1-1}).

Suppose that $m$ is an odd number. Then two equalities given in (ii) can be showed similarly. We omit them here.
\ep

Denote
\begin{eqnarray}\label{e-func-F-f}
g_1(x):=f_1(x,1-x)\quad \mbox{\rm and}\quad g_2(y):=f_2(1-y,y).
\end{eqnarray}
Suppose that $\mathscr{P}=(p_{i_1i_2\cdots i_m})\in \mathbb{T}_{m,2}$ is a symmetric transition probability tensor. Then, it is obvious that (\ref{e-system-1}) has a (unique) solution $(x^*,y^*)^\top\in \Delta$ if and only if $g_1(x)-x=0$ has a (unique) solution $x^*\in [0,1]$, and $1-x^*$ solves $g_2(y)-y=0$. Thus, (\ref{e-system-1}) can be investigated by considering $g_1(x)-x=0$ and $g_2(y)-y=0$. In order to give an appropriate reformulation of the function $g_1(\cdot)$, we need to use the following combinatorial identity.

\begin{lemma}\label{lem-3}
If $n$ is an even number, then
$$
C_n^0+C_n^2+\cdots+C_n^{n}=C_n^1+C_n^3+\cdots+C_n^{n-1}=2^{n-1};
$$
and if $n$ is an odd number, then
$$
C_n^0+C_n^2+\cdots+C_n^{n-1}=C_n^1+C_n^3+\cdots+C_n^n=2^{n-1}.
$$
\end{lemma}

Now, we derive a simple expression of the function $g_1(\cdot)$, which is a key to our discussions later.

\begin{lemma}\label{lem-4}
Suppose that $\mathscr{P}=(p_{i_1i_2\cdots i_m})\in \mathbb{T}_{m,2}$ is a symmetric transition probability tensor, and the function $g_1(\cdot)$ is defined by (\ref{e-func-F-f}). Then, the following results hold.
\begin{itemize}
\item[(i)] If $m$ is an even number, then
$$
g_1(x)=\frac{a-b}{2}(2x-1)^{m-1}+b.
$$
\item[(ii)] If $m$ is an odd number, then
$$
g_1(x)=\frac{a-b}{2}(2x-1)^{m-1}+a.
$$
\end{itemize}
\end{lemma}

\noindent{\bf Proof.} (i) Suppose that $m$ is an even number. In this case, we first show that
\begin{eqnarray}\label{e-lem4-0}
g_1(x)=(a-b)(m-1)\sum_{t\in \{1,3,\ldots, m-2\}}\left\{\frac{2^{m-t-1}}{m-t} C_{m-2}^{t-1}x^{m-t}-\frac{2^{m-t-2}}{m-t-1} C_{m-2}^{t}x^{m-t-1}
\right\}+b. \nonumber\\
\end{eqnarray}
Since $m$ is an even number, it follows from Lemma \ref{lem-1} that
\begin{eqnarray}\label{e-lem4-1}
\begin{array}{rcl}
g_1(x)&=&
ax^{m-1}+bC_{m-1}^1x^{m-2}(1-x)+aC_{m-1}^2x^{m-3}(1-x)^2\\
&&+bC_{m-1}^3x^{m-4}(1-x)^3+\cdots+ax^3(1-x)^{m-4}\\
&&+bC_{m-1}^{m-3}x^2(1-x)^{m-3}+aC_{m-1}^{m-2}x(1-x)^{m-2}+b(1-x)^{m-1}.
\end{array}
\end{eqnarray}
Let
$$
g_1(x)=\sum_{t=1}^{m}\alpha_{m-t}x^{m-t},
$$
then it follows from (\ref{e-lem4-1}) that $\alpha_0=b$, and
\begin{itemize}
\item if $t\in \{1,3,\ldots,m-1\}$, then
$$
\alpha_{m-t}=\sum_{s\in \{0,2,\ldots, m-t-1\}}\left\{C_{m-1}^{t-1+s}C_{t-1+s}^sa-C_{m-1}^{t+s}C_{t+s}^{s+1}b\right\};
$$
\item if $t\in \{2,4,\ldots,m-2\}$, then
$$
\alpha_{m-t}=\sum_{s\in \{2,4,\ldots, m-t-2\}}\left\{C_{m-1}^{t+s}C_{t+s}^sb-C_{m-1}^{t+s+1}C_{t+s+1}^{s+1}a\right\}
+C_{m-1}^{m-1}C_{m-1}^{m-t-1}b.
$$
\end{itemize}
Since
\begin{eqnarray*}
C_{m-1}^{t-1+s}C_{t-1+s}^s
&=& \frac{(m-1)!}{(t-1+s)!(m-t-s)!}\times \frac{(t-1+s)!}{s!(t-1)!}\\
&=& \frac{(m-1)!}{(t-1)!}\times \frac{1}{s!(m-t-s)!}\\
&=& \frac{m-1}{m-t}\times \frac{(m-2)!}{(t-1)!(m-t-1)!}\times \frac{(m-t)!}{s!(m-t-s)!}\\
&=& \frac{m-1}{m-t} C_{m-2}^{t-1}C_{m-t}^s,
\end{eqnarray*}
it follows that for any $t\in \{1,3,\ldots,m-1\}$,
\begin{eqnarray*}
\alpha_{m-t}&=&\frac{m-1}{m-t} C_{m-2}^{t-1}\sum_{s\in \{0,2,\ldots, m-t-1\}}\left\{C_{m-t}^sa-C_{m-t}^{s+1}b\right\}\\
&=& \frac{m-1}{m-t} C_{m-2}^{t-1}[2^{m-t-1}a-2^{m-t-1}b]\\
&=&2^{m-t-1}(a-b)\frac{m-1}{m-t} C_{m-2}^{t-1},
\end{eqnarray*}
where the first equality holds by Lemma \ref{lem-3}; and for any $t\in \{2,4,\ldots,m-2\}$,
$$
\alpha_{m-t}=2^{m-t-1}(b-a)\frac{m-1}{m-t} C_{m-2}^{t-1}.
$$
Thus, (\ref{e-lem4-0}) holds.

Furthermore, it is easy to see that for any $t\in \{1,2,\ldots,m-1\}$,
$$
\frac{m-1}{m-t}C_{m-2}^{t-1}=C_{m-1}^{t-1},
$$
and hence,
\begin{eqnarray*}
g_1(x)&=&(a-b)\sum_{t\in \{1,3,\ldots, m-2\}}\left\{2^{-1} C_{m-1}^{t-1}(2x)^{m-t}-2^{-1} C_{m-1}^{t}(2x)^{m-t-1} \right\}+b\\
&=&\frac{a-b}{2}\sum_{t=1}^{m-2}\left\{C_{m-1}^{t-1}(2x)^{m-t}(-1)^{t-1} \right\}+b\\
&=&\frac{a-b}{2}(2x-1)^{m-1}+b,
\end{eqnarray*}
where the last equality follows from the binomial theorem. Thus, we complete the proof of the result in (i).

(ii) Suppose that $m$ is an odd number. In this case, it follows from Lemma \ref{lem-1} that
\begin{eqnarray}\label{e-lem-4-2}
\begin{array}{rcl}
g_1(x)&:=&
ax^{m-1}+bC_{m-1}^1x^{m-2}(1-x)+aC_{m-1}^2x^{m-3}(1-x)^2\\
&&+bC_{m-1}^3x^{m-4}(1-x)^3+\cdots+ax^4(1-x)^{m-5}+bx^3(1-x)^{m-4}\\
&&+aC_{m-1}^{m-3}x^2(1-x)^{m-3}+bC_{m-1}^{m-2}x(1-x)^{m-2}+a(1-x)^{m-1}.
\end{array}
\end{eqnarray}
Let $g_1(x)=\sum_{t=1}^{m}\beta_{m-t}x^{m-t}$. Then, it follows from (\ref{e-lem-4-2}) that $\beta_0=a$, and
\begin{itemize}
\item if $t\in \{1,3,\ldots,m-2\}$, then
$$
\beta_{m-t}=\sum_{s\in \{0,2,\ldots, m-t-1\}}\left\{C_{m-1}^{t-1+s}C_{t-1+s}^sa-C_{m-1}^{t+s}C_{t+s}^{s+1}b\right\}
+C_{m-1}^{m-1}C_{m-1}^{m-t-1}a;
$$
\item if $t\in \{2,4,\ldots,m-1\}$, then
$$
\beta_{m-t}=\sum_{s\in \{2,4,\ldots, m-t-2\}}\left\{C_{m-1}^{t+s}C_{t+s}^sb-C_{m-1}^{t+s+1}C_{t+s+1}^{s+1}a\right\}.
$$
\end{itemize}
Thus, similar to (i), we can obtain that
$$
g_1(x)=(a-b)(m-1)\sum_{t\in \{1,3,\ldots, m-2\}}\left\{\frac{2^{m-t-1}}{m-t} C_{m-2}^{t-1}x^{m-t}-\frac{2^{m-t-2}}{m-t-1} C_{m-2}^{t}x^{m-t-1}
\right\}+a,
$$
and furthermore,
\begin{eqnarray*}
g_1(x)&=&(a-b)\sum_{t\in \{1,3,\ldots, m-2\}}\left\{2^{-1} C_{m-1}^{t-1}(2x)^{m-t}-2^{-1} C_{m-1}^{t}(2x)^{m-t-1} \right\}+a\\
&=&\frac{a-b}{2}\sum_{t=1}^{m-2}\left\{C_{m-1}^{t-1}(2x)^{m-t}(-1)^{t-1} \right\}+a\\
&=&\frac{a-b}{2}(2x-1)^{m-1}+a,
\end{eqnarray*}
which implies that the result in (ii) holds.
\ep

\begin{lemma}\label{lem-5}
Suppose that $\mathscr{P}=(p_{i_1i_2\cdots i_m})\in \mathbb{T}_{m,2}$ is a symmetric transition probability tensor, and the function $g_1(\cdot)$ is defined by (\ref{e-func-F-f}). Then,
\begin{eqnarray*}
g_1^\prime(x)=(a-b)(m-1)(2x-1)^{m-2}.
\end{eqnarray*}
\end{lemma}

\noindent{\bf Proof.} The desired result follows from Lemma \ref{lem-4} directly.
\ep

The following result is a special case of the one in \cite[Theorem 3.1]{CPZ17}.
\begin{lemma}\label{lem-6}
Suppose that $\mathscr{P}=(p_{i_1i_2\cdots i_m})\in \mathbb{T}_{m,2}$ is a symmetric transition probability tensor. Then, $(\frac 12, \frac 12)^\top$ is a stationary probability vector of $\mathscr{P}$.
\end{lemma}

We now give our main results in this paper.
\begin{theorem}\label{main-thm}
Suppose that $\mathscr{P}=(p_{i_1i_2\cdots i_m})\in \mathbb{T}_{m,2}$ is a symmetric transition probability tensor. If $p_{1\cdots 11}=1$ and $p_{21\cdots 11}=0$, then the corresponding tensor $\mathscr{P}$ is denoted by $\mathscr{P}_1$; and if $p_{1\cdots 11}=0$ and $p_{21\cdots 11}=1$, then the corresponding tensor $\mathscr{P}$ is denoted by $\mathscr{P}_2$. Then, we have the following results.
\begin{itemize}
\item[(i)] The transition probability tensor $\mathscr{P}_1$ has and only has two stationary probability vectors: $(\frac 12,\frac 12)^\top$ and $(0,1)^\top$.
\item[(ii)] The transition probability tensor $\mathscr{P}_2$ has and only has two stationary probability vectors: $(\frac 12,\frac 12)^\top$ and $(1,0)^\top$.
\item[(iii)] If $\mathscr{P}\neq \mathscr{P}_1$ and $\mathscr{P}\neq \mathscr{P}_2$, then $\mathscr{P}$ has a unique stationary probability vector: $(\frac 12, \frac 12)^\top$.
\end{itemize}
\end{theorem}

\noindent{\bf Proof.} Let
\begin{eqnarray}\label{e-g-func}
h(x):=g_1(x)-x.
\end{eqnarray}
Then, it follows from Lemma \ref{lem-5} that
\begin{eqnarray}\label{e-g-func-prime}
h^\prime(x):=g_1^\prime(x)-1=(a-b)(m-1)(2x-1)^{m-2}-1.
\end{eqnarray}
We divide the proof into the following three parts.

\noindent{\bf Part 1}. Suppose that $a=b$. By (\ref{e-g-func-prime}) we have that for any $x\in [0,1]$,
$$
h^\prime(x)=-1<0,
$$
which implies that the function $h(\cdot)$ defined by (\ref{e-g-func}) is strictly decreasing on $[0,1]$. This, together with $h(\frac 12)=0$ by Lemma \ref{lem-1}, implies that $\frac 12$ is the unique solution of $h(x)=0$ on $[0,1]$. Furthermore, by Lemma \ref{lem-1} it follows that $(\frac 12,\frac 12)^\top$ is the unique solution of (\ref{e-system-1}) on $[0,1]$, i.e., $(\frac 12,\frac 12)^\top$ is the unique stationary probability vector of the transition probability tensor $\mathscr{P}$.

\noindent{\bf Part 2}. Suppose that $1\geq a>b\geq 0$. Let $h^\prime(x^*)=0$, then by (\ref{e-g-func-prime}), we have
\begin{eqnarray}\label{e-thm-1}
x^*=\frac 12\left(1+\frac{1}{\sqrt[m-2]{(a-b)(m-1)}}\right).
\end{eqnarray}
We consider the following three cases.
\begin{itemize}
\item[(a)] If $(a-b)(m-1)=1$, then by (\ref{e-thm-1}), we have $x^*=1$. Since $m\geq 3$, it follows that $1>a>b>0$. Since $g_1(1)=a$ by (\ref{e-lem4-1}) and (\ref{e-lem-4-2}), it follows by (\ref{e-g-func}) that $h(1)=a-1<0$, i.e., 1 is not a solution of $h(x)=0$. Moreover, by (\ref{e-g-func-prime}) it follows that the function $h(\cdot)$ is strictly decreasing on $[0,1)$. This, together with $h(\frac 12)=0$, implies that $\frac 12$ is the unique solution of $h(x)=0$ on $[0,1)$. Thus, $\frac 12$ is the unique solution of $h(x)=0$ on $[0,1]$. This, together with Lemma \ref{lem-6}, implies that $(\frac 12,\frac 12)^\top$ is the unique solution of (\ref{e-system-1}) on $[0,1]$, i.e., $(\frac 12,\frac 12)^\top$ is the unique stationary probability vector of the transition probability tensor $\mathscr{P}$.
\item[(b)] If $(a-b)(m-1)<1$, then by (\ref{e-thm-1}), we have $x^*>1$. Thus, by (\ref{e-g-func-prime}) it follows that the function $h(\cdot)$ is strictly decreasing on $[0,1]$. This, together with $h(\frac 12)=0$, implies that $\frac 12$ is the unique solution of $h(x)=0$ on $[0,1]$. Furthermore, by Lemma \ref{lem-6} it follows that $(\frac 12,\frac 12)^\top$ is the unique solution of (\ref{e-system-1}) on $[0,1]$, i.e., $(\frac 12,\frac 12)^\top$ is the unique stationary probability vector of the transition probability tensor $\mathscr{P}$.
\item[(c)] If $(a-b)(m-1)>1$, then by (\ref{e-thm-1}), we have $\frac 12<x^*\leq 1$. Furthermore, by (\ref{e-g-func-prime}) it follows that $h^\prime(x)<0$ when $x\in [0,x^*)$ and $h^\prime(x)>0$ when $x\in (x^*,1]$. On one hand, since the function $h(\cdot)$ is strictly decreasing on $[0,x^*)$ and $h(\frac 12)=0$ with $\frac 12\in [0,x^*)$, it follows that $\frac 12$ is the unique solution of $h(x)=0$ on $[0,x^*)$; and in the meantime, we have $h(x^*)<h(\frac 12)=0$. On the other hand, by (\ref{e-lem4-1}) and (\ref{e-lem-4-2}), we have $g_1(1)=a$, and hence, by (\ref{e-g-func}), $h(1)=a-1\leq 0$. If $a=1$, then since the function $h(\cdot)$ is strictly increasing on $(x^*,1]$ with $h(1)=0$, it follows that $1$ is the unique solution of $h(x)=0$ on $(x^*,1]$; and if $a<1$, then since the function $h(\cdot)$ is strictly increasing on $(x^*,1]$ with $h(1)<0$ and $h(x^*)<0$, it follows that $h(x)=0$ has no solution on $(x^*,1]$. Thus, if $\mathscr{P}=\mathscr{P}_2$, then it has and only has two stationary probability vectors: $(\frac 12,\frac 12)^\top$ and $(1,0)^\top$; otherwise, it has a unique stationary probability vector: $(\frac 12,\frac 12)^\top$.
\end{itemize}

\noindent{\bf Part 3}. Suppose that $1\geq b>a\geq 0$. If $m$ is an even number, then for any $x\in [0,1]$,
$$
h^\prime(x)=(a-b)(m-1)(2x-1)^{m-2}-1<0,
$$
and hence, the function $h(\cdot)$ is strictly decreasing on $[0,1]$. In this case, $(\frac 12,\frac 12)^\top$ is the unique solution of (\ref{e-system-1}) on $[0,1]$, i.e., $(\frac 12,\frac 12)^\top$ is the unique stationary probability vector of the transition probability tensor $\mathscr{P}$.

In the following, we assume that $m$ is an odd number. Let $h^\prime(x^*)=0$, then by (\ref{e-g-func-prime}), we have
\begin{eqnarray}\label{e-thm-2}
x^*=\frac 12\left(1-\frac{1}{\sqrt[m-2]{(b-a)(m-1)}}\right).
\end{eqnarray}
We consider the following three cases.
\begin{itemize}
\item[(a)] If $(b-a)(m-1)=1$, then by (\ref{e-thm-2}), we have $x^*=0$. Since $m\geq 3$, it follows that $1>b>a>0$. Since $m$ is an odd number, we have $g_1(0)=a$ by (\ref{e-lem-4-2}), and hence, it follows by (\ref{e-g-func}) that $h(0)=a-0>0$, i.e., 0 is not a solution of $h(x)=0$. Moreover, by (\ref{e-g-func-prime}) it follows that $h^\prime(x)<0$ on $(0,1]$, which implies that the function $h(\cdot)$ is strictly decreasing on $(0,1]$, and hence, $\frac 12$ is a unique solution of $h(x)=0$ on $(0,1]$. So, $(\frac 12,\frac 12)^\top$ is the unique stationary probability vector of the transition probability tensor $\mathscr{P}$.
\item[(b)] If $(b-a)(m-1)<1$, then by (\ref{e-thm-2}), we have $x^*<0$; and by (\ref{e-g-func-prime}), we have $h^\prime(x)<0$ for any $x\in [0,1]$. Thus, the function $h(\cdot)$ is strictly decreasing on $[0,1]$. This, together with $h(\frac 12)=0$, implies that $\frac 12$ is the unique solution of $h(x)=0$ on $[0,1]$. Thus, $(\frac 12,\frac 12)^\top$ is the unique stationary probability vector of the transition probability tensor $\mathscr{P}$.
\item[(c)] If $(b-a)(m-1)>1$, then it is easy to see from (\ref{e-thm-2}) that $0<x^*<\frac 12$. Furthermore, by (\ref{e-g-func-prime}), we have $h^\prime(x)>0$ when $x\in [0,x^*)$ and $h^\prime(x)<0$ when $x\in (x^*,1]$. On one hand, since the function $h(\cdot)$ is strictly decreasing on $(x^*,1]$ and $h(\frac 12)=0$ with $\frac 12\in (x^*,1]$, it follows that $\frac 12$ is the unique solution of $h(x)=0$ on $(x^*,1]$. Meantime, we have $h(x^*)>h(\frac 12)=0$. On the other hand, since $m$ is an odd number, it follows by (\ref{e-lem-4-2}) that $g_1(0)=a$, and hence, by (\ref{e-g-func}), $h(0)=a\geq 0$. If $a=0$, then $b=1$, and hence, $\mathscr{P}=\mathscr{P}_2$. In this case, $0$ is a solution of $h(x)=0$. Otherwise, since the function $h(\cdot)$ is strictly decreasing on $(0,x^*)$ with $h(0)>0$ and $h(x^*)>0$, it follows that $h(x)=0$ has no solution on $(0,x^*)$. Thus, if $\mathscr{P}=\mathscr{P}_2$, then it has and only has two stationary probability vectors: $(\frac 12,\frac 12)^\top$ and $(1,0)^\top$; otherwise, it has a unique stationary probability vector: $(\frac 12,\frac 12)^\top$.
\end{itemize}
Therefore, by combining {\bf Part 1} with {\bf Part 2} and {\bf Part 3}, we can obtain the desired results.
\ep

By Theorem \ref{main-thm}, we have the following result immediately.
\begin{corollary}\label{coro-1}
Suppose that $\mathscr{P}=(p_{i_1i_2\cdots i_m})\in \mathbb{T}_{m,2}$ is a symmetric transition probability tensor. Then, $\mathscr{P}$ has a unique positive stationary probability vector, which is $(\frac 12,\frac 12)^\top$.
\end{corollary}

Recall that a tensor $\mathscr{P}=(p_{i_1i_2\cdots i_m})\in \mathbb{T}_{m,n}$ is called {\it reducible} if there exists a nonempty proper index subset $I\subset \{1,2,\ldots,n\}$ such that
$$
p_{i_1i_2\cdots i_m}=0,\quad \forall i_1\in I,\;\; \forall i_2,\ldots,i_m\notin I.
$$
If $\mathscr{P}$ is not reducible, then it is called {\it irreducible} \cite{CPZ08}. It is easy to see that both tensors $\mathscr{P}_1$ and $\mathscr{P}_2$ given in Theorem \ref{main-thm} are reducible. Thus, by Theorem \ref{main-thm}, we have the following result immediately.
\begin{corollary}\label{coro-2}
Suppose that $\mathscr{P}=(p_{i_1i_2\cdots i_m})\in \mathbb{T}_{m,2}$ is a symmetric irreducible transition probability tensor. Then, $\mathscr{P}$ has a unique stationary probability vector, which is $(\frac 12,\frac 12)^\top$.
\end{corollary}

When $m=4$, Corollary \ref{coro-2} is just Theorem 3.1 given in \cite{CPZ17}.

\section{Concluding remarks}
\setcounter{equation}{0} \setcounter{assumption}{0}
\setcounter{theorem}{0} \setcounter{proposition}{0}
\setcounter{corollary}{0} \setcounter{lemma}{0}
\setcounter{definition}{0} \setcounter{remark}{0}
\setcounter{algorithm}{0}\setcounter{example}{0}

In this paper, we gave a full characterization on the stationary probability vectors of $m$-order $2$-dimensional symmetric transition probability tensors. In particular, for any integer number $m\geq 3$, any symmetric irreducible transition probability tensor of order $m$ dimension 2 has a unique stationary probability vector. In our analysis, the ``symmetry" of transition probability tensor plays an important role. It is worthy of studying the stationary probability vectors of $m$-order $2$-dimensional transition probability tensors in the absence of symmetry. Moreover, it is also worthy of investigating the uniqueness of the stationary probability vectors of $m$-order $n$-dimensional symmetric transition probability tensors when $m,n\geq 3$.


\begin{thebibliography}{abc99xyz}

\bibitem{BP94}
A. Berman and R. Plemmons, Nonnegative Matrices in the Mathematical Sciences, Classics in
Applied Mathematics, SIAM, 1994.

\bibitem{MZ97}
I. MacDonald and W. Zucchini, Hidden Markov and Other Models for Discrete-valued Time
Series, Chapman \& Hall, London, 1997.

\bibitem{Raf85}
A. Raftery, A model of high-order Markov chains, Journal of the Royal Statistical Society, Series B, 47 (1985): 528-539.

\bibitem{RT94}
A. Raftrey and S. Tavare, Estimation and modelling repeated patterns in high order Markov
chains with the mixture transition distribution model, Applied Statistics, 43 (1994): 179-199.

\bibitem{Wat95}
M. Waterman, Introduction to Computational Biology, Chapman \& Hall, Cambridge, 1995.

\bibitem{CNg06}
W. Ching and M. Ng, Markov Chains: Models, Algorithms and Applications, International
Series on Operations Research and Management Science, Springer, 2006.

\bibitem{CZ13}
K.C. Chang and T. Zhang, On the uniqueness and non-uniqueness of the positive $Z$-eigenvector for transition probability tensors, J. Math. Anal. Appl., 208 (2013): 525-540.

\bibitem{Qi05}
L. Qi, Eigenvalues of a real supersymmetric tensor, J. Symbolic Comput., 40 (2005): 1302-1324.

\bibitem{Lim05}
L.-H. Lim, Singular values and eigenvalues of tensors: a variational approach, in: Proceedings of the IEEE International Workshop on Computational
Advances in Multi-Sensor Adaptive Processing, CAMSAP'05, vol. 1, 2005, pp. 129-132.

\bibitem{NgQZ09}
M. Ng, L. Qi and G. Zhou, Finding the largest eigenvalue of a non-negative tensor, SIAM J. Matrix Anal. Appl., 31 (2009): 1090-1099.



\bibitem{HQ14}
S. Hu and L. Qi, Convergence of a second order Markov chain, Appl. Math. Comput., 241 (2014): 183-192.

\bibitem{LNg14}
W. Li and M. Ng, On the limiting probability distribution of a transition probability tensor, Linear Multilinear Algebra, 62 (2014): 362-385.

\bibitem{LZ16}
C.-K. Li and S. Zhang, Stationary probability vectors of higher-order Markov chains, Linear Algebra Appl., 473 (2016): 114-125.

\bibitem{BB16} H. Bozorgmanesh and M. Hajarian, Convergence of a transition probability tensor
of a higher-order Markov chain to the stationary probability vector, Numer. Linear
Algebra Appl., 23 (2016): 972-988.

\bibitem{CPZ17}
J. Culp, K. Pearson and T. Zhang, On the uniqueness of the $Z_1$-eigenvector of
transition probability tensors, Linear Multilinear Algebra, 65 (2017): 891-896.

\bibitem{CPZ08}
K.C. Chang, K. Pearson, T. Zhang, Perron-Frobenius theorem for nonnegative tensors, Commun. Math. Sci., 6 (2008): 507-520.


\end{thebibliography}
\end{document}